\documentclass{ceb-l}

\usepackage{amssymb}

\newtheorem{theorem}{Theorem}[section]

\theoremstyle{definition}

\theoremstyle{remark}
\newtheorem{remark}[theorem]{Remark}

\numberwithin{equation}{section}

\newcommand{\C}{\mathbf{C}}
\newcommand{\R}{\mathbf{R}}
\newcommand{\sech}{\operatorname{sech}}
\newcommand{\dist}{\operatorname{dist}}
\newcommand{\tr}{\operatorname{tr}}
\newcommand{\eps}{\varepsilon}

\begin{document}


\title{Why are solitons stable?}


\author{Terence Tao}
\address{UCLA Department of Mathematics, Los Angeles, CA 90095-1596.}
\email{tao@@math.ucla.edu}
\thanks{The author is supported by NSF grant CCF-0649473 and a grant from the MacArthur Foundation, and also thanks David Hansen, Robert Miura, Jeff Kimmel, and Jean-Claude Saut for helpful comments and corrections.}


\date{}

\dedicatory{}

\begin{abstract}  The theory of linear dispersive equations predicts that waves should spread out and disperse over time.  However, it is a remarkable phenomenon, observed both in theory and practice, that once nonlinear effects are taken into account, \emph{solitary wave} or \emph{soliton} solutions can be created, which can be stable enough to persist indefinitely.  The construction of such solutions can be relatively straightforward, but the fact that they are \emph{stable} requires some significant amounts of analysis to establish, in part due to symmetries in the equation (such as translation invariance) which create degeneracy in the stability analysis.  The theory is particularly difficult in the \emph{critical} case in which the nonlinearity is at exactly the right power to potentially allow for a self-similar blowup.  In this article we survey some of the highlights of this theory, from the more classical orbital stability analysis of Weinstein and Grillakis-Shatah-Strauss, to the more recent asymptotic stability and blowup analysis of Martel-Merle and Merle-Raphael, as well as current developments in using this theory to rigorously demonstrate controlled blowup for several key equations.
\end{abstract}

\maketitle

\section{Introduction}

In these notes we shall eventually describe recent developments in the stability theory of solitons (or more precisely, \emph{solitary waves}).  Before we discuss solitons, however, we need to first describe the wider context of dispersive equations, and why even the very existence of solitons were initially such a surprising phenomenon.  

In classical physics, it has been realised for centuries that the behaviour of idealised vibrating media (such as waves on string, on the surface of a body of water, or in air), in the absence of friction or other dissipative forces, can be modeled by a number of partial differential equations, known collectively as \emph{dispersive equations}.  Model examples of such equations include the following:

\begin{itemize}
\item The \emph{free wave equation}
$$ u_{tt} - c^2 \Delta u = 0,$$
where $u: \R \times \R^d \to \R$ represents the amplitude $u(t,x)$ of a wave at a point in a spacetime
with $d$ spatial dimensions, $\Delta := \sum_{j=1}^d \frac{\partial^2}{\partial x_j^2}$ is the spatial Laplacian on $\R^d$, $u_{tt}$ is short for $\frac{\partial^2 u}{\partial t^2}$, and $c > 0$ is a fixed constant (which can be rescaled to equal $1$ if one wishes).  This equation models the evolution of waves in a medium which has a fixed speed $c$ of propagation in all directions.
\item The \emph{linear (time-dependent) Schr\"odinger equation}
\begin{equation}\label{schro}
 i \hbar u_t + \frac{\hbar^2}{2m} \Delta u = V u
\end{equation}
where $u: \R \times \R^d \to \C$ is the wave function of a quantum particle, $\hbar, m > 0$ are physical constants (which can be rescaled to equal $1$ if one wishes), and $V: \R^d \to \R$ is a potential function, which we assume to depend only on the spatial variable $x$.  This equation models the evolution of a quantum particle in space in the presence of a classical potential well $V$.
\item The \emph{nonlinear Schr\"odinger (NLS) equation}
\begin{equation}\label{nls}
i u_t + \Delta u = \mu |u|^{p-1} u
\end{equation}
where $p > 1$ is an exponent and $\mu = \pm 1$ is a sign (the case $\mu=+1$ is known as the \emph{defocusing} case, and $\mu=-1$ as the \emph{focusing} case).  This equation can be viewed as a variant of the linear Schr\"odinger equation (with the constants $\hbar$ and $m$ normalised away), in which the potential $V$ now depends in a nonlinear fashion on the solution itself.  This equation no longer has a physical interpretation as the evolution of a quantum particle, but can be derived as a model for quantum media such as Bose-Einstein condensates (see e.g. \cite{spohn}).
\item The \emph{(time-dependent) Airy equation}
\begin{equation}\label{airy}
u_t + u_{xxx} = 0
\end{equation}
where $u: \R \times \R \to \R$ is a scalar function.  This equation can be derived as a very simplified model for propagation of low amplitude water waves in a shallow canal, by starting with the Euler equations, making a number of simplifying assumptions to discard nonlinear terms, and the normalising all constants to equal $1$.  
\item The \emph{Korteweg-de Vries (KdV) equation} \cite{KdV}
\begin{equation}\label{kdv}
 u_t + u_{xxx} + 6 u u_x = 0
\end{equation}
which is a more refined version of the Airy equation in which the first nonlinear term is retained.  The constant $6$ that appears here is not essential, but turns out to be convenient when connecting this equation to the theory of inverse scattering (of which more will be said later).
\item The \emph{generalised Korteweg-de Vries (gKdV) equation}
\begin{equation}\label{gkdv}
 u_t + u_{xxx} + (u^p)_x = 0
\end{equation}
for $p > 1$ an integer; the case $p=2$ is essentially the KdV equation, and the case $p=3$ is known as the \emph{modified Korteweg-de Vries (mKdV) equation}.  The case $p=5$ is particularly interesting due to its \emph{mass-critical} nature, which we will discuss later.
\end{itemize}

For simplicity (to avoid some non-trivial topological phenomena, which are an interesting topic which we have no space to discuss here) we shall only consider equations on flat spacetimes $\R^n$, and only consider solutions which decay to zero at spatial infinity.

The above equations are all evolution equations; if we specify the initial position $u(0,x) = u_0(x)$ of a wave at time $t=0$, we expect these equations to have a unique solution with that initial data\footnote{For the wave equation, which is second-order in time, we also need to specify the initial velocity $\partial_t u(0,x) = u_1(x)$.} for all future times $t > 0$.  Actually, all of the above equations are time-reversible (for instance, if $(t,x) \mapsto u(t,x)$ solves \eqref{kdv}, then $(t,x) \mapsto u(-t,-x)$ also solves \eqref{kdv}) so we also expect the initial data at time $t=0$ to determine the solution at all past times $t < 0$.  (This is in sharp contrast to \emph{dissipative} equations such as the heat equation $u_t = \Delta u$, which are solvable forward in time but are not solvable backwards in time, at least in the category of smooth functions, due to an irreversible loss of energy and information inherent in this equation as time moves forward.)

Solutions to the above equations have two properties, which may seem at first to contradict each other.  The first property is that all of these equations are \emph{conservative}; there exists a \emph{Hamiltonian} $v \mapsto H(v)$, which is a functional that assigns a real number $H(v)$ to any (sufficiently smooth and decaying\footnote{To simplify the exposition, we shall largely ignore the important, but technical, analytic issues of exactly how much regularity and decay one needs in order to justify all the computations and assertions given here.  In practice, one usually first works in the category of \emph{classical solutions} - solutions that are smooth and rapidly decreasing - and then uses rigorous limiting arguments (and in particular, exploiting the \emph{low-regularity well-posedness theory} of these equations) to extend all results to more general classes of solutions, such as solutions in the energy space $H^1(\R^d)$.  See e.g. \cite{cazbooknew}, \cite{kpv:gkdv}, \cite{tao:cbms} for further discussion.}) function $v$ on the spatial domain $\R^d$, such that the Hamiltonian\footnote{Again, with the wave equation, the Hamiltonian depends on the instantaneous velocity $u_t(t)$ of the solution at time $t$ as well as the instantaneous position $u(t)$.} $H(u(t))$ of a (sufficiently smooth and decaying) solution $u(t)$ to the above equation is conserved in time:
$$ H(u(t)) = H(u(0)), \hbox{ or equivalently } \partial_t H(u(t)) = 0.$$
More specifically, the Hamiltonian is given by
$$ H(u(t), u_t(t)) := \int_{\R^d} \frac{1}{2} |u_t(t,x)|^2 + \frac{1}{2} |\nabla_x u(t,x)|^2\ dx$$
for the wave equation,
$$ H(u(t)) := \int_{\R^d} \frac{\hbar^2}{2m} |\nabla_x u(t,x)|^2 + V(x) |u(t,x)|^2\ dx$$
for the linear Schr\"odinger equation,
$$ H(u(t)) := \int_{\R^d} \frac{1}{2} |\nabla_x u(t,x)|^2 + \frac{\mu}{p+1} |u(t,x)|^{p+1}\ dx$$
for the nonlinear Schr\"odinger equation,
$$ H(u(t)) := \int_\R u_x(t,x)^2\ dx$$
for the Airy equation, 
\begin{equation}\label{hut}
H(u(t)) := \int_\R u_x(t,x)^2 - 2 u(t,x)^3\ dx 
\end{equation}
for the Korteweg-de Vries equation, and
\begin{equation}\label{hut-gkdv}
H(u(t)) := \int_\R u_x(t,x)^2 - \frac{2}{p+1} u(t,x)^{p+1}\ dx 
\end{equation}
for the generalised Korteweg-de Vries equation.  In all of these cases, the conservation of the Hamiltonian can be formally verified by computing $\partial_t H(u(t))$ via differentiation under the integral sign, substituting the evolution equation for $u$, and then integrating by parts; we leave this to the reader as an exercise\footnote{One can also formally establish conservation of the Hamiltonian by interpreting each of the above dispersive equations in turn as an infinite-dimensional Hamiltonian system, but we will not adopt this (important) perspective here; see \cite{Kuksin} for further discussion.}.

Actually, the Hamiltonian is not the only conserved quantity available for these equations; each of these equations also enjoy a number of symmetries (e.g. translation invariance), which (by Noether's theorem) leads to a number of important additional conserved quantities, which we will discuss later.  It is often helpful to interpret the Hamiltonian as describing the total \emph{energy} of the wave.

The conservative nature of these equations means that even for very late times $t$, the state $u(t)$ of solution is still\footnote{This is assuming that the solution exists all the way up to this time $t$, which can be a difficult task to  establish rigorously, especially if the initial data was rough.  Again, we suppress these important technical issues for simplicity.} ``similar'' to the initial state $u(0)$, in the sense that they have the same energy.  This may lead one to conclude that solutions to these evolution equations should evolve in a fairly static, or fairly periodic manner; after all, this is what happens to solutions to finite-dimensional systems of ordinary differential equations which have a conserved energy $H$ which is \emph{coercive} in the sense that the \emph{energy surfaces} $\{ H = \hbox{const}\}$ are always bounded.

However, this intuition turns out to not be correct in the realm of dispersive equations, even though such equations can be thought of as infinite-dimensional systems of ODE with a conserved energy, and even though this energy usually exhibits coercive properties.  This is ultimately because of a second property of all of these equations, namely \emph{dispersion}.  Informally, dispersion means that different components of a solution $u$ to any of these equations travel at different velocities, with the velocity of each component determined by the frequency.  As a consequence, even though the state of solution at late times has the same energy as the initial state, the different components of the solution are often so dispersed that the solution at late times tends to have much smaller amplitudes than at early times\footnote{This phenomenon may seem to be inconsistent with time reversal symmetry.  However, this dispersive effect only occurs when the initial data is spatially localised; dispersion sends localised high-amplitude states to broadly dispersed, low-amplitude states, but (by time reversal) can also have the reverse effect.}.  Thus, for instance, it is perfectly possible for the solution $u(t)$ to go to zero in $L^\infty(\R^d)$ norm as $t \to \pm \infty$, even as its energy stays constant (and non-zero).  The ability to go to zero as measured in one norm, while staying bounded away from zero in another, is a feature of systems with infinitely many degrees of freedom, which is not present when considering systems of ODE with only boundedly many degrees of freedom.

One can see this dispersive effect in a number of ways.  One (somewhat informal) way is to analyse \emph{plane wave} solutions
\begin{equation}\label{utx}
 u(t,x) = A e^{i t \tau + x \cdot \xi}
\end{equation}
for some non-zero amplitude $A$, some temporal frequency $\tau \in \R$, and some spatial frequency $\xi \in \R^d$.  For instance, for the Airy equation \eqref{airy}, one easily verifies that \eqref{utx} solves\footnote{Strictly speaking, one needs to allow solutions $u$ to \eqref{airy} to be complex-valued here rather than real-valued, but this is of course a minor change.} \eqref{airy} exactly when $\tau = \xi^3$; this equation is known as the \emph{dispersion relation} for the Airy equation.  If we rewrite the right-hand side of \eqref{utx} in this case as $A e^{i \xi (x - (- \xi^2) t)}$, this asserts that a plane wave solution to \eqref{airy} has a \emph{phase velocity} $-\xi^2$ which is the negative square of its spatial frequency $\xi$.  Thus we see that for this equation, higher frequency plane waves have a much faster phase velocity than lower frequency ones, and the velocity is always in a leftward direction.  Similar analyses can be carried out for the other equations given above, though in those equations involving a nonlinearity or a potential, one has to restrict attention to small amplitude or high frequency solutions so that one can (non-rigorously) neglect the effect of these terms.  For instance, for the Schr\"odinger equation \eqref{schro} (at least with $V=0$) one has the dispersion relation 
\begin{equation}\label{disp}
\tau = - \frac{\hbar^2}{2m} |\xi|^2.
\end{equation}

However, as is well known in physics, the phase velocity does not determine the speed of propagation of information in a system; that quality is instead controlled by the \emph{group velocity}, which typically is slightly different from the phase velocity.  To explain this quantity, let us modify the ansatz \eqref{utx} by allowing the amplitude $A$ to vary slowly in space, and propagate in time at some velocity $v \in \R^d$.  More precisely, we consider solutions of the form
$$ u(t,x) \approx A( \eps (x - vt ) ) e^{it\tau + x \cdot \xi} $$
where $A$ is a smooth function and $\eps > 0$ is a small parameter, and we shall be vague about what the symbol ``$\approx$'' means.  If we have $\tau = \xi^3$ as before, then a short (and slightly non-rigorous) computation shows that
$$ u_t + u_{xxx} \approx \eps (v + 3 \xi^2) A'( \eps(x-vt) ) e^{it\tau + x \cdot \xi} + O(\eps^2).$$
Thus we see that in order for $u$ to (approximately) solve \eqref{airy} up to errors of $O(\eps^2)$, the group velocity $v$ must be equal to $-3\xi^2$, which is three times the phase velocity $-\xi^2$.  Thus, at a qualitative level at least, we still have the same predicted behaviour as before; all frequencies propagate leftward, and higher frequencies propagate faster than lower ones.  In particular we expect localised high-amplitude states, which can be viewed (via the Fourier inversion formula) as linear superpositions of plane waves of many different frequencies, to disperse leftwards over time into broader, lower-amplitude states (but still with the same energy as the original state, of course).

One can perform similar analyses for other equations.  For instance, for the linear Schr\"odinger equation, and assuming either high frequencies or small potential, one expects waves to propagate at a velocity proportional to their frequency, according to \emph{de Broglie's law} $mv = \hbar \xi$; similarly for the nonlinear Schr\"odinger equation when one assumes either high frequencies or small amplitude.  In contrast, for wave equation, this analysis suggests that waves of (non-zero) frequency $\xi$ should propagate at velocities $\frac{\xi}{|\xi|} c$; thus the propagation speed $c$ is constant but the propagation direction $\frac{\xi}{|\xi|}$ varies with frequency, leading to a weak dispersive effect.  For more general dispersive equations, the group velocity can be read off of the dispersion relation $\tau = \tau(\xi)$ by the formula $v = -\nabla_\xi \tau$ (whereas in contrast, the phase velocity is $-\frac{\xi}{|\xi|^2} \tau$).

In the case of the Schr\"odinger equation with $V=0$, one can see the dispersive effort more directly by using the explicit solution formula
$$ u(t,x) = \frac{1}{(2\pi \hbar t/m)^{d/2}} \int_{\R^d} e^{i\frac{m|x-y|^2}{2\hbar^2 t}} u_0(y)\ dy$$
for $t \neq 0$ and all sufficiently smooth and decaying initial data $u_0$.  Indeed, we immediately conclude from this formula (formally, at least) that if $u_0$ is absolutely integrable, then $\|u(t)\|_{L^\infty(\R^d)}$ decays at the rate $O( |t|^{-d/2} )$.

In the case of linear equations such as the Airy equation \eqref{airy}, there is a similar explicit formula (involving the Airy function $\operatorname{Ai}(x)$ instead of complex exponentials), but one can avoid the use of special functions by instead proceeding using the Fourier transform and the principle of stationary phase (see e.g. \cite{stein:large}).  Indeed, by starting with the Fourier inversion formula
$$ u_0(x) = \int_\R \hat u_0(\xi) e^{i x \xi}\ d\xi$$
where $\hat u_0(\xi) := \frac{1}{2\pi} \int_\R u_0(x) e^{-ix\xi}\ dx$ is the Fourier transform of $u_0$, and noting as before that $e^{it\xi^3 + ix\xi}$ is the solution of the Airy equation with initial data $e^{ix\xi}$, we see from the principle of superposition (and ignoring issues of interchanging derivatives and integrals, etc.) that the solution $u$ is given by the formula
\begin{equation}\label{utx-airy}
u(t,x) = \int_\R \hat u_0(x) e^{it\xi^3 + i x\xi}\ d\xi.
\end{equation}
If $u_0$ is a Schwartz function (infinitely smooth, with all derivaives decreasing faster than any polynomial), then its Fourier transform is also Schwartz and thus slowly varying.  On the other hand, as $t$ increases, the phase $e^{it\xi^3 + i \xi}$ oscillates more and more rapidly (for non-zero $\xi$), and so we expect an increasing amount of cancellation in the integral in \eqref{utx-airy}, leading to decay of $u$ as $t \to \infty$.  This intuition can be formalised using the methods of stationary phase (which can be viewed as advanced applications of the undergraduate calculus tools of integration by parts and changes of variable), and can for instance be used to show that $\|u(t)\|_{L^\infty(\R)}$ decays at a rate $O(t^{-1/3})$ in general.

This technique of representing a solution as a superposition of plane waves also works (with a twist) for the linear Schr\"odinger equation \eqref{schro} in the presence of a potential $V$, provided that the potential is sufficiently smooth and decaying.  The basic idea is to replace the plane waves \eqref{utx} by \emph{distorted plane waves} $\Phi(\tau,x) e^{it\tau}$, where (in order to solve \eqref{schro}) $\Phi$ has to solve the \emph{time-independent Schr\"odinger equation}
\begin{equation}\label{phiv}
 - \tau \hbar \Phi + \frac{\hbar^2}{2m} \Delta \Phi = V \Phi,
 \end{equation}
and then to try to represent solutions $u$ to \eqref{schro} as superpositions
$$ u(t,x) = \int a(\tau) \Phi(\tau,x) e^{it\tau}\ d\tau$$
where we are being intentionally vague as to what the range of integration is.  If we restrict attention to negative values of $\tau$, then it turns out (by use of scattering theory) that we can construct distorted plane waves $\Phi(\tau,x)$ which asymptotically resemble the standard plane waves $e^{i x \cdot \xi}$ as $|x| \to \infty$, where $\xi$ is a frequency obeying the dispersion relation
\eqref{disp}.  If $u$ is composed entirely of these waves, then one has a similar dispersive behaviour to the free Schr\"odinger equation (for instance, under suitable regularity and decay hypotheses on $V$ and $u_0$, $\|u(t)\|_{L^\infty(\R^d)}$ will continue to decay like $O(t^{-d/2})$).  In such cases we say that $u$ is in a \emph{radiating state}.  In many important cases (such as when the potential $V$ is non-negative, or is small in certain function space norms), all states (with suitable regularity and decay hypotheses) are radiating states.  However, when $V$ is large and allowed to be negative, it is also possible\footnote{When $V$ does not decay rapidly, then there can also be some intermediate states involving the singular continuous spectrum of the Schr\"odinger operator $-\frac{\hbar^2}{2m} \Delta + V$, which disperse over time slower than the radiating states but faster than the bound states.  One can also occasionally have \emph{resonances} corresponding to the boundary case $\tau=0$, which exhibit somewhat similar behaviour.  For simplicity of exposition, we will not discuss these (important) phenomena.}
 to contain \emph{bound states}, in which $\tau$ is positive, and the distorted plane wave $\Phi(\tau,x)$ is replaced by an \emph{eigenfunction} $\Phi$, which continues to solve the equation
\eqref{phiv}, but now $\Phi$ decays exponentially to zero as $|x| \to \infty$, instead of oscillating like a plane wave as before.  (Informally, this is because once $\tau$ is positive, the dispersion relation \eqref{disp} is forcing $\xi$ to be imaginary rather than real.)  In particular, $\Phi$ lies in $L^2(\R^d)$, and so $-\tau$ becomes an eigenvalue of the Schr\"odinger operator\footnote{This operator $H$ is related to the Hamiltonian $H(u)$ discussed earlier by the formula $H(u) = \langle Hu, u\rangle$, where $\langle u,v \rangle := \int_{\R^d} u \overline{v}$ is the usual inner product on $L^2(\R^d)$.} $H := -\frac{\hbar^2}{2m} \Delta + V$.  Because multiplication $V$ is a compact operator relative to 
$-\frac{\hbar^2}{2m} \Delta$, standard spectral theory shows that the set of eigenvalues $-\tau$ is discrete (except possibly at the origin $-\tau=0$).  Note that it is necessary for $V$ to take on negative values in order to obtain negative eigenvalues, since otherwise the operator $H$ is positive semi-definite.

If $u_0$ consists of a superposition of one or more of these eigenfunctions, e.g.
$$ u_0 = \sum_k c_k \Phi(\tau_k,x)$$
where $-\tau_k$ ranges over finitely many of the eigenvalues of $H$, then we formally have
$$ u(t) = \sum_k c_k e^{it\tau_k} \Phi(\tau_k,x),$$
and so we see that $u(t)$ oscillates in time but does not disperse in space.  In this case we say that $u$ is a \emph{bound state}.  Indeed, the evolution is instead \emph{almost periodic}, in the sense that $\liminf_{t \to \infty} \|u(t)-u_0\|_{L^2(\R^d)} = 0$, or equivalently that the orbit $\{ u(t): t \in \R \}$ is a precompact subset of $L^2(\R^d)$.

By further application of spectral theory, one can show that an arbitrary state $u_0$ (in, say, $L^2(\R^d)$) can be decomposed as the orthogonal sum of a radiating state, which disperses as $t \to \infty$, and a bound state, which evolves in an almost periodic manner.  Indeed this decomposition corresponds to the decomposition of the spectral measure of $H$ into absolutely continuous and pure point components.

\section{Solitons}

We have seen how solutions to linear dispersive equations either disperse completely as $t \to \infty$, or else (in the presence of an external potential) decompose into a superposition of a radiative state that disperses to zero, plus a bound state that exhibits phase oscillation but is otherwise stationary.

In everyday physical experience with water waves, we of course see that such waves disperse to zero over time; once a rock is thrown into a pond, for instance, the amplitude of the resulting waves diminish over time.  However, one does see in nature water waves which refuse to disperse for astonishingly long periods of time, instead moving at a constant speed without change in shape.  Such \emph{solitary waves} or \emph{solitons}\footnote{Strictly speaking, a wave which is localised and maintains its form for long periods of time is merely a \emph{soilitary wave}.  A soliton is a solitary wave with the additional property that solitons and other radiation can pass through it without destroying its form.  The question of understanding the collision between two solitary waves for non-integrable equations is still poorly understood despite some recent progress, so we shall focus here instead on the more perturbative question of how solitary waves interact with small amounts of radiation.  In the literature on that subject, it is then customary to refer to the solitary wave as a soliton, though this is technically an abuse of notation.}  were first observed by John Scott Russell, who followed such a wave in a shallow canal on horseback for over a mile, and then reproduced such a traveling wave (which he called a ``wave of translation'') in a wave tank.

This phenomenon was first explained mathematically by Korteweg and de Vries \cite{KdV} in 1895, using the equation \eqref{kdv} that now bears their name (although this equation was first proposed as a model for shallow wave propagation by Boussinesq a few decades earlier).  Indeed, if one considers traveling wave solutions to \eqref{kdv} of the form
$$ u(t,x) = f(x-ct)$$
for some velocity $c$, then this will be a solution to \eqref{kdv} as long as $f$ solves the ODE
$$ - c f' + f''' + 6 f f' = 0.$$
If we assume that $f$ decays at infinity, then we can integrate this third-order ODE to obtain a second-order ODE
$$ - cf + f'' + 3 f^2 = 0.$$
For $c > 0$, this ODE admits the localised explicit solutions $f(x) = c Q(c^{1/2} (x-x_0))$ for any $x_0 \in \R$, where $Q$ is the explicit Schwartz function $Q(x) := \frac{1}{2} \sech^2(\frac{x}{2})$.  For $c \leq 0$, one can show that there are no localised solutions other than the trivial solution $f \equiv 0$.  Thus we obtain a family of explicit solitary wave solutions 
\begin{equation}\label{soliton}
u(t,x) = c Q(c^{1/2}( x - ct - x_0 ) )
\end{equation}
to the KdV equation; the parameter $c$ thus controls the speed, amplitude, and width of the wave, while $x_0$ determines the initial location.  

Interestingly, all the solutions \eqref{soliton} move to the \emph{right}, while radiating states move to the left.  This phenomenon is somewhat analogous to the situation with the linear Schr\"odinger equation, in which the temporal frequency $\tau$ (which is somewhat like the propagation speed $c$ in KdV) was negative for radiating states and positive for bound states.  Similar solitary wave solutions can also be found for gKdV and NLS equations, though in higher dimensions $d > 1$ one cannot hope to obtain such explicit formulae for these solutions, and instead one needs to use more modern PDE tools, such as calculus of variations and other elliptic theory methods, in order to build such solutions; see e.g. \cite{blions}, \cite{blions2}, \cite{gidas}.  There are also larger and more oscillatory ``excited'' solitary wave solutions which, unlike the ``ground state'' solitary wave solutions described above, exhibit changes of sign, but we will not discuss them here.

Early numerical analyses of the KdV equation\cite{zab}, \cite{fermi} revealed that these soliton solutions \eqref{soliton} were remarkably stable.  Firstly, if one perturbed a soliton by adding a small amount of noise, then the noise would soon radiate away from the soliton, leaving the soliton largely unchanged (other than some slight perturbation in the $c$ and $x_0$ parameters); these phenomena are described mathematically by results on the \emph{orbital stability} and \emph{asymptotic stability} of solitons, of which more will be said later.  This is perhaps unsurprising, given that solitons move rightwards and radiation moves leftwards, but one has to bear in mind that equations such as \eqref{kdv} are not linear, and in particular one cannot obviously superimpose a soliton and a radiative state to create a new solution to the KdV equation.

What was even more surprising was what happened if one considered collisions between two solitons, for instance imagining initial data of the form
$$ u(0,x) = c_1 Q(c_1^{1/2}( x - x_1 ) ) + c_2 Q(c_2^{1/2}( x - x_2 ) )$$
with $0 < c_2 < c_1$ and $x_1$ far to the left of $x_2$; thus initially we have a larger, fast-moving soliton to the left of a shallower, slow-moving soliton.  If the KdV equation were linear, the solution would now take the form
$$ u(t,x) = c_1 Q(c_1^{1/2}( x - c_1 t - x_1 ) ) + c_2 Q(c_2^{1/2}( x - c_2 t - x_2 ) )$$
and so the faster solitons would simply overtake the slower one, with no interaction between the two.  At the other extreme, with a strongly nonlinear equation, one could imagine all sorts of scenarios when two solitons collide, for instance that they scatter into radiation or into many smaller solitons, combine into a large soliton, and so forth.  However, the KdV equation exhibits an interesting intermediate behaviour: the solitons do interact nonlinearly with each other during collision, but then emerge from that collision almost unchanged, except that the solitons have been shifted slightly by their collision.  In other words, for very late times $t$, the solution approximately takes the form
$$ u(t,x) \approx c_1 Q(c_1^{1/2}( x - c_1 t - x_1 - \theta_1 ) ) + c_2 Q(c_2^{1/2}( x - c_2 t - x_2 - \theta_2) )$$
for some additional shift parameters $\theta_1, \theta_2 \in \R$.

More generally, if one starts with \emph{arbitrary} (but smooth and decaying) initial data, what usually happens (numerically, at least) with evolutions of equations such as \eqref{kdv} is that some non-linear (and chaotic-seeming) behaviour happens for a while, but eventually most of the solution radiates away to infinity and a finite number of solitons emerge, moving away from each other at different rates.  Quite remarkably, this behaviour can in fact be justified rigorously for the KdV equation and a handful of other equations (such as the NLS equation in the cubic one-dimensional case $d=1, p=3$) due to the inverse scattering method, which we shall discuss shortly, although even in those cases, there are some exotic solutions, such as ``breather'' solutions, which occasionally arise and which do not evolve to a superposition of solitons and radiation, but instead exhibit periodic or almost periodic behaviour in time.  Nevertheless, it is widely believed (and supported by extensive numerics) that for many other dispersive equations (roughly speaking, those equations whose nonlinearity is not strong enough to cause finite time blowup, and more precisely for the \emph{subcritical} equations), solutions with ``generic'' initial data should eventually resolve into a finite number of solitons, moving at different speeds, plus a radiative term which goes to zero.  This (rather vaguely defined) conjecture goes by the name of the \emph{soliton resolution conjecture}.  Except for those few equations which admit exact solutions (for instance, by inverse scattering methods), the conjecture remains unsolved in general, in part because we have very few tools available that can say anything meaningful about \emph{generic} data in a certain class (e.g. with some function norm bounds) without also being applicable to \emph{all} data in that class; thus the presence of a few exotic solutions that do not resolve into solitons and radiation seems to prevent us from tackling all the other cases.  Nevertheless, there are certain important regimes in which we do have a good understanding.  One of these is the perturbative regime near a single soliton, in which the initial state $u_0$ is close to that of a soliton such as \eqref{soliton}; this case will be the main topic of our discussion.  More recently, results have begun to emerge on \emph{multisoliton} states, in which the solution is close to the superposition of many widely separated solitons, and even more recently still, there has been some results on the collision between a very fast narrow soliton and a very slow broad one.  However, it seems that truly non-perturbative regimes, such as the collisions between two solitons of comparable size, remain beyond the reach of current tools (perhaps requiring a new advance in our understanding of dynamical systems in general).  (See \cite{tao:soliton}, \cite{soffer-icm}, \cite{tao:soliton2} for further discussion.)

\section{The inverse scattering approach}

We now briefly mention the technique of \emph{inverse scattering}, which is a non-perturbative approach which allows one to control the evolution of solutions to completely integrable equations such as \eqref{kdv}.  (Similar techniques apply to one-dimensional cubic NLS, see e.g. \cite{akns}, \cite{sa}, \cite{manakov}, \cite{noveksenov}.) This is a vast subject, which can be viewed from many different algebraic and geometric perspectives; we shall content ourselves with describing the approach based on \emph{Lax pairs} \cite{lax}, which has the advantage of simplicity, provided that one is willing to accept a rather miraculous algebraic identity.

The identity in question is as follows.  Suppose that $u$ solves the KdV equation \eqref{kdv}.  As always we assume enough smoothness and decay to justify the computations that follow.  For every time $t$, we consider the time-dependent differential operators $L(t)$, $P(t)$ acting on functions on the real line $\R$, defined by
\begin{align*}
L(t) &:= \frac{d^2}{dx^2} + u(t) \\
P(t) &:= \frac{d^3}{dx^3} + \frac{3}{4} ( \frac{d}{dx} u(t) + u(t) \frac{d}{dx})
\end{align*}
where we view $u(t)$ as a multiplication operator, $f \mapsto u(t) f$.  One can view $P(t)$ as a truncated (non-commutative) Taylor expansion of $L(t)^{3/2}$.  In view with this interpretation, it is perhaps not so surprising that $L(t)$ and $P(t)$ ``almost commute''; the commutator $[P(t),L(t)] := P(t) L(t) - L(t) P(t)$ of the third order operator $P(t)$ and the second order operator $L(t)$ would normally be expected to be fourth order, but in fact things collapse to just be zeroth order.  Indeed, after some computation, one eventually obtains
$$ [P(t),L(t)] = \frac{1}{4}( u_{xxx}(t) + 6 u(t) u_x(t) ).$$
In particular, if we substitute in \eqref{kdv}, we obtain the remarkable \emph{Lax pair equation}
\begin{equation}\label{lp}
\frac{d}{dt} L(t) = 4 [P(t), L(t)].
\end{equation}
If we non-rigorously treat the operators $L(t)$, $P(t)$ as if they were matrices, we can interpret this equation as follows.  Using the Newton approximation
$$ L(t+dt) \approx L(t) + dt \frac{d}{dt} L(t); \quad \exp( \pm P(t) dt ) \approx (1 \pm P(t) dt )$$
for infinitesimal $dt$, we see from \eqref{lp} that
\begin{equation}\label{ltdt}
L(t+dt) \approx \exp( 4 P(t) dt ) L(t) \exp( -4 P(t) dt ).
\end{equation}
This informal analysis suggests that $L(t+dt)$ is a conjugate of $L(t)$, and so on iterating this we expect $L(t)$ to be a conjugate of $L(0)$.  In particular, the \emph{spectrum} of $L(t)$ should be time-invariant.  Since $L(t)$ is determined by $u(t)$, this leads to a rich source of invariants for $u(t)$.

The above analysis can be made more rigorous.  For instance, one can show that the traces\footnote{Actually, to avoid divergences we need to consider normalised traces, such as $\tr( e^{sL(t)} - e^{s \frac{d^2}{dx^2}} )$.} $\tr(e^{sL(t)})$ of heat kernels are independent of $t$ for any fixed $s > 0$; expanding those traces in powers of $s$ one can recover an infinite number of conservation laws, which includes the conservation of the Hamiltonian \eqref{hut} as a special case.  We will not pursue this approach further here, and refer the reader to \cite{hitchin}.  Another way to proceed is to consider solutions to the generalised eigenfunction equation
\begin{equation}\label{lphi}
L(t) \phi(t,x) = \tau \phi(t,x)
\end{equation}
for some $\tau \in \R$ and some smooth function $\phi(t,x) = \phi(t,x;\tau)$ (not necessarily decaying at infinity).  If the equation
\eqref{lphi} holds for a single time $t$ (e.g. $t=0$), and if $\phi$ then evolves by the equation
\begin{equation}\label{phit}
\phi_t(t,x) = 4 P(t) \phi(t,x)
\end{equation}
for all $t$, one can verify (formally, at least) that \eqref{lphi} persists for all $t$, by differentiating \eqref{lphi} in time and substituting in \eqref{lp} and \eqref{phit}.  (The astute reader will note that these manipulations are equivalent to those used to produce \eqref{ltdt}.  

This now suggests an strategy to solve the KdV equation exactly from an arbitrary choice of initial data $u(0) = u_0$.

\begin{enumerate}
\item Use the initial data $u_0$ to form the operator $L(0)$, and then locate the generalised eigenfunctions $\phi(0,x;\lambda)$ for each choice of spectral parameter $\tau$.
\item Evolve each generalised eigenfunction $\phi$ in time by the equation \eqref{phit}.
\item Use the generalised eigenfunctions $\phi(t,x;\tau)$ to recover $L(t)$ and $u(t)$.
\end{enumerate}

This strategy looks very difficult to execute, because the operator $P(t)$ itself depends on $u(t)$, and so \eqref{phit} cannot be solved exactly without knowing what $u(t)$ is - which is exactly what we are trying to find in the first place!  But we can break this circularity by only seeking to solve \eqref{phit} \emph{at spatial infinity} $x = \pm \infty$.  Indeed, if $u(t)$ is decaying, and $\tau = -\xi^2$ for some real number $\xi$, then we see that solutions $\phi(t,x)$ to \eqref{lphi} must take the form
$$ \phi(t,x) \approx a_\pm(\xi;t) e^{i\xi x} + b_\pm(\xi;t) e^{-i\xi x}$$
as $x \to \pm \infty$, for some quantities $a_\pm(\xi;t), b_\pm(\xi;t) \in \C$, which we shall refer to as the \emph{scattering data} of $L(t)$.  (One can normalise, say, $a_-(0) = 1$ and $b_-(\xi;0) = 0$, and focus primarily on $a_+(\xi;t)$ and $b_+(\xi;t)$, if desired.)  Applying \eqref{phit} and using the decay of $u(t)$ once again, we are then led (formally, at least) to the asymptotic equations
$$ \partial_t a_\pm(\xi;t) = 4 (i\xi)^3 a_\pm(\xi;t); \quad \partial_t b_\pm(\xi;t) = 4 (-i\xi)^3 b_\pm(\xi;t)$$
which can be explicitly solved\footnote{Note the resemblance of the phases here to those in \eqref{utx-airy}.  This is not a co-incidence, and indeed the scattering and inverse scattering transforms can be viewed as nonlinear versions of the Fourier and inverse Fourier transform.};
\begin{equation}\label{at-evolve}
a_\pm(\xi;t) = e^{-4i\xi^3 t} a_\pm(\xi;0); \quad b_\pm(\xi;t) = e^{-4i\xi^3 t} b_\pm(\xi;0).
\end{equation}
This only handles the case of negative energies $\lambda < 0$.  For positive energies, say $\lambda = +\xi^2$ for some $\xi > 0$, the situation is somewhat similar; in this case, we have a discrete set of $\xi$ for which we have a decaying solution $\phi(t,x)$, with $\phi(t,x) \approx c_\pm(\xi;t) e^{\mp \xi x}$ for $x \to \pm \infty$, where 
\begin{equation}\label{ct-evolve}
c_\pm(\xi;t) = e^{\mp 4\xi^3 t} c_\pm(\xi;0).
\end{equation}

This suggests a revised strategy to solve the KdV equation exactly:

\begin{enumerate}
\item Use the initial data $u_0$ to form the operator $L(0)$, and then locate the scattering data $a_\pm(\xi;0)$, $b_\pm(\xi;0)$, $c_\pm(\xi;0)$.
\item Evolve the scattering data by the equations \eqref{at-evolve}, \eqref{ct-evolve}.
\item Use the scattering data at $t$ to recover $L(t)$ and $u(t)$.
\end{enumerate}

The main difficulty in this strategy is now the third step, in which one needs to solve the \emph{inverse scattering problem} to recover $u(t)$ from the scattering data.  This is a vast and interesting topic in its own right, and involves complex-analytic problems such as the Riemann-Hilbert problem; we will not discuss it further here (but see e.g. \cite{ggkm}, \cite{akns}, \cite{dt}).  Suffice to say, though, that after some work, it is possible to execute the above strategy for sufficiently smooth and decaying initial data $u$ to obtain what is essentially\footnote{The solution is not quite expressible as a closed-form integral, as in \eqref{utx-airy}, but can be built out of solving a number of ordinary differential-integral equations (such as the Gelfand-Levitan-Marchenko equation, see e.g. \cite{dt}), which turns out to suffice for the purposes of analysing the asymptotic behaviour of the solution.} an explicit formula for $u$.

The relationship of all this to solitons is as follows.  Recall from our discussion of the linear Schr\"odinger equation \eqref{schro} that the operator $L(0) = \frac{d^2}{dx^2} + u_0$ is going to have radiating states (or absolutely continuous spectrum) corresponding to negative energies $\tau = -\xi^2 < 0$, and a discrete set of positive eigenfunctions corresponding to positive energies $\tau = +\xi^2 > 0$.  Generically, the eigenvalues are simple\footnote{Repeated eigenvalues lead to more complicated behaviour, including breather solutions and logarithmically divergent solitons.}.  In that case, it turns out that the inverse scattering procedure relates each eigenvalue $+\xi^2$ of $L(0)$ to a soliton present inside $u_0$; the value of $\xi$ determines the scaling parameter $c$ of the soliton, and the scattering data $c_\pm(\xi;0)$ determines (in a slightly complicated fashion, depending on the rest of the spectrum) the location of the solitons.  The remaining scattering data $a_\pm(\xi;0)$, $b_\pm(\xi;0)$ determines the radiative portion of the solution.  As the solution evolves, the spectrum stays constant, but the data $a_\pm, b_\pm, c_\pm$ changes in a controlled manner; this is what causes the solitons to move and the radiation to scatter.  It turns out that the exact location of each soliton depends to some extent on the relative sizes of the constants $c_\pm$, which are growing or decaying exponentially at differing rates; it is because of this that as one soliton overtakes another, the location of each soliton gets shifted slightly.  

\section{The analytic approach}

The inverse scattering method gives extremely powerful and precise information on very general (and in particular, non-perturbative) solutions to equations such as the Korteweg-de Vries equation.  However, it does not seem to be directly applicable to more general equations, such as the gKdV equation \eqref{gkdv} for\footnote{The modified KdV equation with $p=3$ turns out to also be completely integrable; in fact, it can even be transformed directly into the KdV equation by a simple operation known as the \emph{Miura transform}\cite{miura}, which we will not discuss further here.} $p \neq 2,3$.  For instance, no reasonable Lax pair formulation exists for these equations.  We now turn to more analytic techniques, which are less sensitive to the fine algebraic structure of the equation, although they still do rely very heavily on conservation laws and their relatives, such as monotonicity formulae.

We shall mostly restrict attention to the gKdV equation \eqref{gkdv}.  We have already identified one conserved quantity for this equation, namely the energy \eqref{hut-gkdv}.  Another such conserved quantity is the mass
$$ M(u(t)) := \int_\R u(t,x)^2\ dx.$$
Together, the mass and energy can (in some cases) control the $H^1$ norm
$$ \| u(t) \|_{H^1_x(\R)}^2 := \int_\R u(t,x)^2 + u_x(t,x)^2\ dx.$$
Indeed, if we are in the \emph{mass-subcritical} case $p < 5$, then the \emph{Gagliardo-Nirenberg inequality}
\begin{equation}\label{gn}
\int_\R v^{p+1} \leq C(p) (\int_\R v^2)^{\frac{p+3}{4}} (\int_\R v_x^2)^{\frac{p-1}{4}},
\end{equation}
valid for any $v$ with suitable decay and regularity, gives us the \emph{a priori} bound
\begin{equation}\label{muh}
 \| u(t) \|_{H^1_x(\R)} \leq C( M(u(t)), H(u(t)) ) = C( M(u_0), H(u_0) )
\end{equation}
for some quantity $C( M(u_0), H(u_0) )$ depending on the initial mass or energy.  The condition $p<5$ is necessary to ensure that the exponent $\frac{p-1}{4}$ in \eqref{gn} is strictly less than $1$.  This condition can also be deduced from scaling heuristics, by investigating how the mass and energy transform under the scale invariance
\begin{equation}\label{scaling}
u(t,x) \mapsto \lambda^{-\frac{2}{p-1}} u(\frac{t}{\lambda^3}, \frac{x}{\lambda}).
\end{equation}

It is possible to use the \emph{a priori} bound \eqref{muh}, combined with the Picard iteration method for constructing solutions, and some moderately advanced estimates from harmonic analysis, to show that the equation \eqref{gkdv} in the mass-subcritical case admits a unique global smooth solutions from arbitrary smooth, decaying data; see \cite{kpv:gkdv}.  Thus there is no problem with existence, uniqueness, or regularity when it comes to these equations; the only remaining analytic issue (albeit a difficult one) is to understand the asymptotic behaviour of solutions.

The above analysis for $p < 5$ is valid no matter how large the mass $M(u) = M(u_0)$ of the solution is.  If we then turn to the \emph{mass-critical} case $p=5$, the situation changes; the \emph{a priori} bound \eqref{muh} is now only valid when the mass $M(u_0)$ is sufficiently small.  In fact, by using the sharp Gagliardo-Nirenberg inequality of Weinstein\cite{Weinstein}, one can be more precise as follows.  Given any $p > 1$, the equation \eqref{gkdv} admits a family of soliton
(or solitary wave) solutions similar to \eqref{soliton}, namely
\begin{equation}\label{soliton-gkdv}
u(t,x) = c^{1/(p-1)} Q( c^{1/2}( x - x_0 - ct ) )
\end{equation}
where 
$$ Q(x) := \left( \frac{p+1}{2 \cosh^2( \frac{p-1}{2} x )} \right)^{1/(p-1)}$$
is a positive, smooth, rapidly decreasing solution to the ODE 
\begin{equation}\label{qxx}
Q_xx + Q^p = Q.
\end{equation}
(In fact, up to translation, $Q$ is the only such solution; see \cite{coff}, \cite{kwong}.)  In particular, we have the \emph{standard soliton} solution to gKdV
\begin{equation}\label{standard}
u(t,x) = Q(x-t);
\end{equation}
all other solitons differ from the standard soliton only by the scaling \eqref{scaling} and the translation invariance.

In the mass-critical case $p=5$, all of the solitons have the same mass, namely
$$ M(u) = M(Q) = \int_\R Q^2.$$
The Gagliardo-Nirenberg inequality of Weinstein can then be used to show that one has the a priori bound \eqref{muh} in the $p=5$ case as long as one only considers solutions with mass strictly less than that of the ground state, and so long as the solution exists; see \cite{Weinstein}.  The latter caveat is a substantial one; it is conjectured that one has global existence of solutions for gKdV of the $p=5$ from smooth decaying initial data whenever the mass is less than that of the ground state, but this conjecture is still open.  (However, global existence is known if the mass is sufficiently small, by a perturbative argument based on the contraction mapping principle and some harmonic analysis estimates; see \cite{kpv:gkdv}.)  An important recent work of Martel and Merle \cite{liouville}, \cite{merle}, \cite{mm-annals}, \cite{mm-jams}, though, shows in this case that singularities can form for data arbitrarily close to the ground state (but of slightly larger mass); we discuss this result in more detail in Section \ref{crit}.

In the \emph{mass-supercritical case} $p > 5$ the situation is very unclear, due to the lack of good \emph{a priori} etimates in this case; it is likely that singularities do form in this case for large initial data, but this has not been rigorously established.  It is however known that solitons are unstable in this setting.

We return now to the mass-subcritical case $p<5$, in which global existence and regularity are assured, and now consider the problem of \emph{stability} of one of the soliton solutions \eqref{soliton-gkdv}.  By taking advantage of the scaling and translation invariance in the problem, we can reduce matters to considering the stability of the standard soliton $Q(x-t)$.  For instance, if we are given a solution $u$ which is close to this soliton at time zero, for instance in the sense that
\begin{equation}\label{uc0}
\| u(0) - Q \|_{H^1(\R)} \leq \sigma
\end{equation}
for some sufficiently small $\sigma$, is this enough to guarantee that $u$ stays close to \eqref{soliton} for much later times, thus
\begin{equation}\label{uc1}
\| u(t) - Q(\cdot - t ) ) \|_{H^1(\R)} \leq \sigma'
\end{equation}
for some other small $\sigma'$ depending on $\sigma$, and for large times $t$?  (For small times $t$, the local well-posedness theory allows one to obtain bounds of this form, but with $\sigma'$ replaced by $C \sigma \exp(Ct)$ for some constant $C$ depending only on $p$.)  One can phrase this question in other norms than the $H^1$ norm, of course, but this norm turns out to be rather natural due to its connection with the Hamiltonian (which we have already seen in \eqref{muh}).

This type of \emph{absolute stability} of the soliton is too strong a property to hold, basically because it is not compatible with the scale invariance \eqref{scaling}.  Indeed, consider the soliton solution \eqref{soliton-gkdv} with $x_0=0$ and $c = 1 + O(\sigma)$ very close to $1$. Then \eqref{uc0} holds, but \eqref{uc1} will fail for sufficiently large times $t$, because $u$ has most of its mass (and $H^1$ norm) near $ct$, whereas \eqref{soliton} has most of its mass near $ct$.  The point is that by rescaling the soliton very slightly, one can adjust the speed of that soliton, which will eventually over time cause the perturbed soliton to diverge from the original soliton.  Note that this conclusion has nothing to do with the $H^1$ norm, and would work for basically any reasonable function space norm.

However, even though this perturbed soliton is far from the original soliton at late times $t$, it is still close to a \emph{translation} of that original soliton (by $ct-t$).  Equivalently, if we define the \emph{ground state curve} 
$$\Sigma = \{ Q(\cdot - x_0): x_0 \in \R \} \subset H^1(\R)$$
consisting of all translates of the ground state\footnote{More suggestively, one should think of $\Sigma$ as the space of all possible soliton states whose conserved statistics (in particular, mass and energy) agree with that of the ground state.  In the case of the sub-critical NLS equation \eqref{nls}, $\Sigma$ then becomes a cylinder, formed by considering the action of both translation $Q(x) \mapsto Q(x-x_0)$ and phase rotation $Q(x) \mapsto e^{i\theta} Q(x)$ on $Q$.  In the critical case, the dimension of $\Sigma$ increases due to the additional symmetry of scale invariance, as we shall shortly see.} $Q$, then we see that $u(t)$ stays close to $\Sigma$ for all $t$.  To put it another way, while $u(t)$ does not stay close to $Q(\cdot-t)$ for each time $t$, the \emph{orbit} $\{ u(t): t \in \R \}$ stays close to the orbit $\{ Q(\cdot-t): t \in \R \} = \Sigma$.  Indeed, this is a general phenomenon:

\begin{theorem}[Orbital stability of sub-critical gKdV]\label{orb}\cite{benjamin}, \cite{bona}, \cite{wein:modulate}  Let $1 < p < 5$.  If $u_0 \in H^1(\R)$ is such that $\dist_{H^1}(u_0,\Sigma)$ is sufficiently small (say less than $\sigma$ for some small constant $\sigma > 0$), and $u$ is the solution to \eqref{gkdv} with initial data $u_0$, then we have
$$ \dist_{H^1}(u(t),\Sigma) \lesssim \dist_{H^1}(u_0, \Sigma)$$
for all $t$.  Here we use $X \lesssim Y$ or $X = O(Y)$ to denote the estimate $|X| \leq CY$ for some $C$ that depends only on $p$, and $X \sim Y$ as shorthand for $X \lesssim Y \lesssim X$.
\end{theorem}

This theorem is proven by a variant of the classical \emph{Lyapunov functional} method for establishing absolute stability.  Let us briefly recall how that method works.  Suppose we were able to find a functional $u \mapsto L(u)$ on $H^1$ with the following properties:

\begin{enumerate}
\item If $u$ is an $H^1$ solution to \eqref{gkdv}, then $L(u(t))$ is non-increasing in $t$.  
\item $Q$ is a local minimiser of $L$, thus $L(u) - L(Q) \geq 0$ for all $u$ sufficiently close to $Q$ in $H^1$.
\item Furthermore, the minimum is non-degenerate in the sense that $L(u)-L(Q) \sim \| u-Q\|_{H^1}^2$ for all $u$ sufficiently close to $Q$ in $H^1$.  
\end{enumerate}

These three facts would then easily imply that $Q$ is absolutely stable.  Indeed, if $u_0$ is close to $Q$, then $L(u_0)$ is close to (but not smaller than) $L(Q)$, which implies that $L(u(t))$ is also close to but not smaller than $L(Q)$ for all $t > 0$, which implies (by a continuity argument) that $u(t)$ is close to $Q$ for all $t > 0$.  (The case $t<0$ can then be handled by time reversal symmetry $u(t,x) \mapsto u(-t,-x)$, or equivalently by considering $L(u(-\cdot))$ instead of $L(u(\cdot))$.)

We already saw, though, that $Q$ is not absolutely stable, and so such a Lyapunov functional cannot exist.  However, we can still hope to obtain a modified Lyapunov functional $u \mapsto L(u)$ which implies orbital stability instead of absolute stability.  More precisely, we require $L$ to be such that

\begin{enumerate}
\item If $u$ is an $H^1$ solution to \eqref{gkdv}, then $L(u(t))$ is non-increasing in $t$.  
\item $L(u) = L(Q)$ for all $u \in \Sigma$.
\item $\Sigma$ is a local minimiser of $L$, thus $L(u) - L(Q) \geq 0$ for all $u$ sufficiently close to $\Sigma$ in $H^1$.
\item Furthermore, the minimum is non-degenerate in the sense that $L(u)-L(Q) \sim \dist_{H^1}(u,\Sigma)^2$ for all $u$ sufficiently close to $\Sigma$ in $H^1$.  
\end{enumerate}

It is not hard to see that this would imply Theorem \ref{orb}. The task then reduces to locating the functional $L$ with the stated properties.  From properties 2,4 it seems reasonable to look for an $L$ which is translation invariant.  From property 1 and time reversal symmetry it seems reasonable to look for an $L$ which is conserved, such as a combination of the mass $M(u)$ and the energy $H(u)$.  It also has to be a functional for which $Q$ is a local minimum, thus \eqref{qxx} should essentially be the Euler-Lagrange equation for $L$.  With all these heuristics, one is soon led to the candidate
$$ L(u) := H(u) + M(u).$$
It is then not hard to verify most of the required properties for $L$, especially if we \emph{define} $Q$ to be the minimiser of $L$.  The one tricky thing is to show the strict non-degeneracy $L(u) - L(Q) \gtrsim \dist_{H^1}(u,\Sigma)^2$ when $u$ is close to $\Sigma$.  One difficulty here is the translation invariance of the estimate; if we do not break this symmetry, then we are forced to only use translation-invariant methods to establish the estimate, which greatly reduces the range of tools available.  Hence we shall \emph{break} the symmetry by decomposing
$$ u = Q(\cdot - x_0) + \eps$$
for some small function\footnote{Note here that $\eps = \eps(t,x)$ is denoting a \emph{function} rather than a number!  This notation is traditional in the literature.} $\eps \in H^1(\R)$ and some $x_0 \in \R$.  There are a number of ways we can choose the parameter $x_0$.  The most obvious approach is to pick $Q(\cdot-x_0)$ to be the translated ground state which is closest to $u$ in $H^1$ norm, thus minimising $\|\eps\|_{H^1}$.  By elementary calculus, this allows us to obtain the orthogonality condition
$$ \langle \eps, Q'(\cdot-x_0) \rangle_{H^1(\R)} = 0$$
where $\langle u,v\rangle_{H^1(\R)} := \int_\R uv + u_x v_x$ is the standard inner product on $H^1$.  Other choices of $x_0$ will lead to a slightly different orthogonality condition; some orthogonality conditions are more suitable for some applications than others, but we will not explore this technical issue further here.  

We can then break (or ``spend'') the translation invariance by normalising $x_0$ to be zero, thus $u = Q+\eps$ and $\eps$ is orthogonal to $Q'$.  (Note that $Q'$ represents the infinitesimal action of the translation group at $Q$.)  Now, since $Q$ is a minimiser of $L$ we have (formally, at least) the Taylor expansion
$$ L(Q+\eps) = L(Q) + \frac{1}{2} L''(Q)(\eps,\eps) + O( \eps^3 )$$
where $L''(Q): H^1(\R) \times H^1(\R) \to \R$ is some explicit positive semi-definite symmetric bilinear form.  The task is then to show that
$$ L''(Q)(\eps,\eps) \gtrsim \|\eps\|_{H^1}^2$$
when $\eps$ is orthogonal to $Q'$.  (An orthogonality condition of this sort is necessary; since $L$ is translation invariant, we easily verify that $L''(Q)$ must annihilate $Q'$.)  This is a \emph{spectral gap} condition on $L''(Q)$, which can be viewed as a positive-definite self-adjoint operator, and can be established by spectral methods; the key ingredient needed is a \emph{uniqueness} result that asserts that $Q$ and its translates are the \emph{only} minimisers of $L$.  Details can be found \cite{wein:modulate}\footnote{If one replaces the power nonlinearity in \eqref{gkdv} with a more general nonlinearity, then this positive definiteness of $L''(Q)$ can fail.  In that case one can in fact show that the solitons are not orbitally stable.}.

The theory of orbital stability for very general dispersive models has now been extensively developed, see e.g. \cite{gss}, \cite{bss}.

Another way to state the above results is that if a global solution $u$ starts off close to the ground state curve $\Sigma$, then at later times one has the decomposition
\begin{equation}\label{que}
u(t,x) = Q(x-x(t)) + \eps(t, x-x(t))
\end{equation}
for some function $x: \R \to \R$ (which tracks the position of the soliton component of $u$) and some error term $\eps$, which is small in $H^1$.  We have the freedom to impose one (non-degenerate) orthogonality condition of our choice on $\eps$, such as $\langle \eps, Q' \rangle_{H^1} = 0$, by choosing $x(t)$ appropriately.

The question then arises as to what happens to the error $\eps$ over time, or to the position $x(t)$.  We return to the model example of the rescaled soliton \eqref{soliton-gkdv}.  In this case we can take $x(t) = x_0 + ct$ and $\eps(t,x) = c^{1/(p-1)} Q( c^{1/2} x ) - Q(x)$ (this $\eps$ does not quite obey the above orthogonality condition, but this will not concern us).  Thus we see in this case that $\eps$ does not disperse to zero in any sense.  However, we can hope to ``quotient out'' this scaling and obtain a decomposition \eqref{que} which has a better error term $\eps$.  Indeed, if we replace the ground state curve $\Sigma$ by the \emph{ground state surface}
\begin{equation}\label{sigmap}
 \Sigma' := \{ c^{1/(p-1)} Q( c^{1/2} (\cdot - x_0) ): c > 0, x_0 \in \R \} \supset \Sigma
 \end{equation}
and approximate $u(t)$ by an element of $\Sigma'$, we can obtain a more refined decomposition
\begin{equation}\label{utx-asym}
u(t,x) = R(t,x) + \eps(t,x)
\end{equation}
where $R$ is a soliton-like state
\begin{equation}\label{qcc}
R(t,x) := c(t)^{1/(p-1)} Q( c(t)^{1/2}( x - x(t) ) ) 
\end{equation}
for some \emph{modulation functions} $x: \R \to \R$ and $c: \R \to \R^+$ (with $c$ always close to $1$), and $\eps$ is small in $H^1$ (one can take $\|\eps\|_{H^1}=O(\sigma)$ if $\dist_{H^1}(u_0,\Sigma) \leq \sigma$).  The point is now that by enlarging the dimension of the approximating surface from $1$ to $2$, the error $\eps$ is now allowed to enjoy \emph{two} orthogonality conditions rather than just one.  There are again several choices of which orthogonality conditions to pick (anything which is suitably ``transverse'' to $\Sigma'$ will do); a typical set of choices is 
\begin{equation}\label{rtx}
 \int_\R R(t,x) \eps(t,x)\ dx = \int_\R (x-x(t)) R(t,x) \eps(t,x)\ dx = 0.
\end{equation}

Now we can hope that with such a refined decomposition that the error $\eps$ will disperse, especially in the neighbourhood of $x(t)$, so that in the vicinity of the solition $R(t,x)$, the solution $u$ converges locally to $R(t,x)$.  Let us informally say that the soliton is \emph{asymptotically stable} if we have a result of this form.  Such stability results can be obtained for the completely integrable cases $p=2,3$ by using the inverse scattering methods of the previous section.  For general sub-critical $p$, such results were first obtained by Pego and Weinstein \cite{pego} and by Mizumachi \cite{mizu}, for perturbations of the ground state which were strongly localised (e.g. assuming exponential decay at infinity).  More recently, Martel and Merle \cite{martel}, \cite{martel2} were able to consider more general perturbations which were only assumed to be small in the energy norm $H^1$; this generalisation is important for the purposes of understanding how a soliton will collide with another shallow broad soliton, which may have small energy but will not have strong localisation properties.  In particular, they showed

\begin{theorem}[Asymptotic stability for sub-critical gKdV]\label{asymt}  Let the notation and assumptions be as in Theorem \ref{orb}.  Then we have a decomposition of the form \eqref{utx-asym}, with $c(t) = c_+$ constant and close to $1$, $x(t)$ differentiable with 
\begin{equation}\label{xpt}
\lim_{t \to +\infty} x'(t) = c_+,
\end{equation}
and the error term $\eps$ obeying the local decay $\| \eps(t) \|_{H^1( x > \beta t)} \to 0$ as $t \to +\infty$ for any $\beta > 0$.
\end{theorem}

Roughly speaking, this asserts that as $t \to \infty$, the solution resolves into a soliton moving to the right at an asymptotically constant speed $c_+$, plus an error term which is radiating to the left; this is of course consistent with the soliton resolution conjecture.  (In the case $p=4$, and with the additional scale-invariant assumption that $u-Q$ is small in $\dot H^{1/6}$, a refinement of this result was given by the author, asserting that $\eps$ in fact converges asymptotically to a solution of the Airy equation \eqref{airy}.)

The estimate \eqref{xpt} implies the asymptotic $x(t) = c_+ t + o(t)$.  It is not entirely clear what the nature of the $o(t)$ error is; one might naively expect to obtain a refined asymptotic of the form $x(t) = c_+ t + x_+ + o(1)$, but it turns out that by inverse scattering methods one can give an example in the $p=2$ case in which one has the asymptotic $x(t) = t + \kappa \sqrt{\log t} + o(\sqrt{\log t})$ for some $\kappa > 0$.

We now sketch the ideas used to prove Theorem \ref{asymt}.  The first step is to pass from \eqref{gkdv}, which is an equation describing the dynamics of $u$, to equations describing the dynamics of $\eps$, $x(t)$, and $c(t)$.  This can be done by algebraic manipulations\footnote{As always, we ignore the analytic issues of how to justify all the formal computations in the case when $u$ is low regularity; this can be done by standard (and boring) regularisation and limiting arguments.}.  Indeed, if one substitutes \eqref{utx-asym}, \eqref{qcc} into \eqref{gkdv}, one eventually obtains the equation
\begin{equation}\label{eps}
\eps_t + \eps_{xxx} + (p R^{p-1} \eps)_x = F_{x(t),c(t),x'(t),c'(t)} + N(\eps, R)
\end{equation}
for $\eps$ where the forcing term $F_{x(t),c(t),x'(t),c'(t)}$ (caused by changes in the modulation parameters) is the explicit smooth function 
\begin{equation}\label{force}
F_{x(t),c(t),x'(t),c'(t)} := -\frac{c'(t)}{c(t)} (\frac{2R}{p-1} + (x-x(t)) R_x) + (x'(t)-c(t)) R_x
\end{equation}
and $N(\eps,R)$ (caused by self-interactions of the radiation term $\eps$) is the nonlinearity
$$ N(\eps,R) := ( (R+\eps)^p - R^p - p R^{p-1} \eps )_x.$$
As for the evolution of $x(t)$ and $c(t)$, one can differentiate \eqref{rtx} to obtain a $2 \times 2$ linear system of equations (known as the \emph{modulation equations}) expressing the evolution $x'(t)$ and $c'(t)$ of the modulation parameters in terms of various integrals involving $R, \eps$ and its derivatives.  The exact form of these modulation equations is not important for our purposes; the only thing which matters is the type of control that one gets on $x'(t)$ and $c'(t)$.  By comparison with the soliton solutions \eqref{soliton-gkdv} one expects $x'(t)$ to be close to $1$, and $c'(t)$ to be close to $0$.  For most choices of orthogonality conditions, the degree of this closeness will only be linear in $\eps$.  But if one uses the specific orthogonality conditions \eqref{rtx}, it turns out that there are particular cancellations which allow the error here to be \emph{quadratic} in $\eps$, at least as regards the variation of the scale parameter $c(t)$.  Indeed, one can show after some computation (exploiting the exponential decay of $R$ and its derivatives away from $x(t)$) that
\begin{equation}\label{cpt}
|c'(t)| + |x'(t)-1|^2 \lesssim \int_\R \eps^2(t,x) e^{-|x-x(t)|}\ dx.
\end{equation}
This is a rather strong estimate; it asserts that the error term $\eps$ only has a linear influence on the velocity, and a quadratic influence on the change in scale, and only when a significant portion of the mass of $\eps$ is stationed near the soliton.  These bounds are particularly useful in controlling the size of the forcing term $F_{x(t),c(t),x'(t),c'(t)}$.

The right-hand side of \eqref{eps} now consists primarily of terms which behave quadratically or higher in $\eps$.  This raises the hope that one can use perturbation theory to approximate the evolution here by that of the linearised equation $\eps_t + \eps_{xxx} + (p R^{p-1} \eps)_x = 0$.  (There is still one term, namely the drift term $(x'(t)-c(t)) R_x$in \eqref{force}, in the right-hand side which exhibits linear behaviour, but this term only causes a translation in $\eps$ and is thus be manageable.)  To do this, we need\footnote{If we do not exploit this fact, then our control on the dispersive effects of the linearised equation is too weak; we can only hope to obtain decay of $O(t^{-1/3})$ on $\eps$ at best, which is insufficient to allow us to neglect the quadratic nonlinearity terms.} to somehow exploit the fact that the linearised equation is trying to propagate $\eps$ to the left, while the soliton is moving to the right.

One strategy to achieve this is via an understanding of the linearised equation $\eps_t + \eps_{xxx} + (p R^{p-1} \eps)_x = 0$.  After some rescaling one can replace $R$ with $Q$ here.  If one works in suitably weighted spaces, then one can use spectral theory to obtain good decay properties for this evolution which can be used to neglect the nonlinear terms and recover asymptotic stability in the case of rapidly decreasing perturbations.  However, this approach seems to become very delicate in the case of perturbations in the energy space.

One particularly elegant way to achieve control on the error $\eps$ for perturbations in the energy space is via \emph{virial identities}, as carried out in \cite{martel2}.  Let us motivate these identities in the simple model case of the Airy equation \eqref{airy}.  This equation has a conserved mass, indeed one quickly computes using \eqref{airy} and integration by parts that
$$ \partial_t \int_\R u^2 = \int_\R 2 u u_t = - \int_\R 2u u_{xxx} = 0.$$
To affirm the intuition that the mass of $u$ should be propagating leftward, let us now introduce the virial quantity $\int_\R x u^2$, which one can think of as the mean position of $u$.  We compute
\begin{align*}
\partial_t \int_\R x u^2 &= \int_\R 2 x u u_t \\
&= - \int_\R 2x u u_{xxx} \\
&= \int_\R 2 u u_{xx} + 2x u_x u_{xx} \\
&= - \int_\R 3 u_x^2. 
\end{align*}
In particular, we see that $\int_\R x u^2$ is a decreasing function of time, which is a quantitative realisation of the intuition of leftwards propagation.  If we instead replace $x$ by $x-x(t)$, we get even faster decay:
$$ \partial_t \int_\R (x-x(t)) u^2 = - \int_\R 3 u_x^2 - x'(t) \int_\R u^2.$$
In particular, if $x'(t) \gtrsim 1$ (which is the situation we are in above), we have
$$ \partial_t \int_\R (x-x(t)) u^2 \leq - c \|u\|_{H^1(\R)}^2$$
for some $c > 0$.

It turns out that one can do the same sort of thing for \eqref{eps}.  Indeed, one can show (after lengthy computations) that (formally, at least) we have
\begin{equation}\label{epsor}
 \partial_t \int_\R (x-x(t)) \eps^2 \leq -c \|\eps \|_{H^1(\R)}^2
\end{equation}
for some $c > 0$.  This estimate strongly suggests that $\eps$ will move to the left of the soliton over time.  Unfortunately, this ``global'' virial identity cannot be used directly in the above analysis, because the integral on the left-hand side may be divergent due to lack of spatial decay on $\eps$.  However, this can be rectified by the usual trick of localising the weight $x-x(t)$.  Indeed, one can show that for sufficiently large $A > 1$, we have the ``local'' virial identity
$$ \partial_t \int_\R \Psi_A(x-x(t)) \eps^2 \leq -c \int_\R (\eps^2 + \eps_x^2) e^{-|x-x(t)|/A}$$
for some bounded increasing function $\Psi(x-x(t))$ which equals $x-x(t)$ for $|x-x(t)| \leq A$, and is of magnitude $O(A)$ throughout.  Thus the quantity $\int_\R \Psi_A(x-x(t)) \eps^2$ is monotone decreasing, while also being controlled by $A$ times the mass.  If we then integrate this in time we obtain an important spacetime bound
\begin{equation}\label{epsx}
\int_\R \int_\R (\eps^2 + \eps_x^2) e^{-|x-x(t)|/A}\ dx dt \lesssim A \sigma.
\end{equation}
This is the first indication of dispersion away from the soliton; it asserts that the radiation term $\eps$ cannot linger near the soliton $x(t)$ for extended periods of time.

This estimate, combined with \eqref{cpt}, is already enough to demonstrate convergence of the scale parameter $c(t)$ to an asymptotic limit $c_+$.  It also shows that the forcing term in \eqref{eps} decays quite quickly in time; in particular, the quadratic nature of the nonlinearity shows that it decays integrably in time, with the exception of the
drift term $(x'(t)-c(t)) R_x$ which can be dealt with by hand.  Because of this, it is possible to use energy estimates to conclude the full strength of Theorem \ref{asymt}.

We now describe an alternate approach to asymptotic stability, also due to Martel and Merle \cite{martel}, which uses a more sophisticated and general strategy which has since shown to be useful for many other equations, including critical equations.  The basic strategy is to use the \emph{compactness-and-contradiction method}, which we informally summarise as follows.

\begin{enumerate}
\item Suppose we wish to show some asymptotic property of a solution $u(t)$ as $t \to +\infty$.  We assume for contradiction that this property does not occur.
\item By using weak compactness, we then extract a sequence $t_n$ of times going to infinity in which (suitably normalised versions) of the state $u(t_n)$ are weakly convergent in some sense, but which violate the property in some quantitative manner.  In particular, (suitably normalised versions of) $u(t+t_n)$ should converge weakly to some asymptotic solution $u_\infty(t)$ of the original equation (now defined for all times $t \in \R$), which continues to violate the desired property.
\item By using the dispersive properties of the equation, show that the asymptotic solution $u_\infty(t)$ obeys some \emph{strong} compactness properties, or equivalently that the evolution $t \mapsto u_\infty(t)$ is \emph{almost periodic} in some strong topology.   (At this point $u_\infty$ is behaving somewhat like the dispersive analogue of a solution to an elliptic PDE or variational problem.)
\item Using more dispersive properties of the equation, upgrade the strong compactness to obtain further regularity and decay of the solution.  (This step is roughly analogous to the exploitation of \emph{elliptic regularity} in the theory of elliptic PDE.)
\item Establish a \emph{Liouville theorem} or \emph{rigidity theorem}, that the only solutions close to solitons which exhibit strong compactness, regularity, and decay properties are the solitons itself.  This is the most difficult step, and often requires full use of the conservation laws and monotonicity formulae of the equation.  (This is analogous to Liouville theorems in elliptic PDE, the most famous of which is the assertion that the only bounded holomorphic or harmonic functions on $\C$ or $\R^d$ are the constants.)
\item We conclude that $u_\infty$ is itself a soliton, which we then combine with the fact that it violates the required property to obtain a contradiction.
\end{enumerate}

The compactness-and-contradiction method is extremely powerful in analysing many nonlinear parabolic and dispersive equations, for instance a variant of this method for Ricci flow also plays a crucial role in Perelman's recent proof of the Poincar\'e conjecture.  Another variant of this method is also very useful in establishing large data global well-posedness results for critical equations, though we will not discuss this topic further here.
The one drawback of the method is that, by being indirect and relying so strongly on compactness methods, it does not easily provide any sort of quantitative bound in its conclusions, in contrast to the previous arguments used to prove Theorem \ref{asymt}, which were direct and easily provide explicit bounds.

Let us now sketch how this method is applied to give a new proof of Theorem \ref{asymt}.  Actually we will just prove the slightly weaker claim that the translated radiation terms $\eps(t, x-x(t))$ converges weakly in $H^1(\R)$ to zero as $t \to +\infty$; note this weak convergence implies for instance that $\eps(t,x-x(t))$ converges locally uniformly to zero, and so the radiation term eventually vacates the neighbourhood of the soliton.  One can upgrade this convergence to obtain results closer in strength in Theorem \ref{asymt}, but we will not do so here.

To prove this weak convergence claim, we use the compactness-and-contradiction method.  Suppose for contradiction that $\eps(t, x-x(t))$ does not converge weakly to $H^1(\R)$ as $t \to +\infty$.  Since $\eps$ is bounded in $H^1$, weak compactness then shows that there exists a sequence of times $t_n \to \infty$ such that $\eps(t_n, x-x(t_n))$ converges weakly in $H^1$ to some non-trivial limit $\eps_\infty(0, x)$; one can also assume that $c(t_n)$ converges to some limit $c_+$.  Due to some weak continuity properties of the gKdV flow (which can be proven by harmonic analysis methods) one can then show that $u(t+t_n,x+x(t_n))$ converges weakly (and locally in time) to some limiting solution $u_\infty(t,x) = R_\infty(t,x) + \eps_\infty(t,x)$, where $R_\infty$ and $\eps_\infty$ obey similar estimates to $R$ and $\eps$, and $R_\infty$ is defined using some modulation parameters $c_\infty(t)$ and $x_\infty(t)$.

The normalised radiation terms $\eps_\infty(t, x-x(t))$ stay bounded in $H^1$.  By the Rellich compactness theorem, this means that they are locally precompact in $L^2$, i.e. their restriction to any compact spatial interval $I$ lies in a compact subset of $L^2(I)$.  We now assert that these terms are in fact \emph{globally} precompact in $L^2$.  This is equivalent to asserting that for any $\delta > 0$, we must have some radius $R$ such that we have very little mass on the left,
\begin{equation}\label{xl}
 \int_{x < x(t) - R} |\eps_\infty(t,x)|^2\ dx < \delta
\end{equation}
and very little mass on the right
\begin{equation}\label{xr}
\int_{x > x(t) + R} |\eps_\infty(t,x)|^2\ dx < \delta.
\end{equation}

We briefly sketch why one would expect these claims to be true.  Suppose that \eqref{xr} failed, then a non-zero portion of the mass of $\eps_\infty$ at some time would be far to the right of the soliton.  Returning to the original solution, we see that there exist arbitrarily large times $t$ for which a significant portion of the mass of $\eps(t)$ is to the right of $x(t)$.  Now we evolve backwards in time, back to time $0$.  Away from the soliton, mass has a tendency to move leftwards as one goes forwards in time, and thus rightwards as one goes backwards in time.  One can make this precise (by using crude forms of the local virial identity alluded to before), and conclude that at time $0$, a significant portion of the mass of $\eps(0)$ is to the right of $x(t)$.  But $t$ can be arbitrarily large, and so $x(t)$ can be arbitrarily large also (recall that $x'(t)$ stays close to $1$).  This contradicts the monotone convergence theorem, and so \eqref{xr}.

The proof of \eqref{xl} is similar; if \eqref{xl} failed, then at some point a non-zero portion of the mass of $\eps_\infty$ lies far to the left of the soliton, and thus we have strictly less than $M(\eps_\infty)$ of the mass of $\eps_\infty$ near or to the right of the soliton.  By the above discussion, we see that we have strictly less than $M(\eps_\infty)$ of the mass of $\eps(t)$ near or to the right of $x(t)$ for a sequence of arbitrarily large times $t$.  But by using local virial-type identities to control the propagation of mass, this loss of mass to the left is irreversible, and in fact we have strictly less than $M(\eps_\infty)$ of the mass of $\eps(t)$ near or to the right of $x(t)$ for \emph{all} sufficiently large times $t$.  But then it is not possible for $\eps(t+t_n,x+x(t_n))$ to converge weakly to $\eps_\infty(t,x)$, a contradiction.  

It turns out that one can upgrade the bounds \eqref{xl}, \eqref{xr} significantly, to obtain a pointwise uniform exponential decay estimate of the form
\begin{equation}\label{epst}
 |\eps_\infty(t,x)| \lesssim \|\eps_\infty(t)\|_{H^1} e^{-c |x-x(t)|}
\end{equation}
for some $c > 0$.   This is established by a long-time analysis of \eqref{eps} and exploits the fact that the fundamental solution to the Airy equation \eqref{airy} decays exponentially fast in the rightwards direction.  This uniformity estimate is crucial in what follows.

It then remains to establish the \emph{Liouville theorem} that if $\eps_\infty$ is a sufficiently small (in $H^1$) solution to a nonlinear equation \eqref{eps} with $\eps_\infty(t,x-x(t))$ compact in $L^2$, then $\eps_\infty$ must vanish.  To prove this, we first use another compactness-and-contradiction argument in order to eliminate the nonlinear terms in \eqref{eps}. If the claim failed, then we could find a sequence of solutions $\eps_n$ to \eqref{eps} which
converged to zero in $H^1$ norm as $n \to \infty$, and were each compact in $L^2$, but were non-zero.  Normalising each $\eps_n$ by its $H^1$ norm, we see from from the uniform estimate \eqref{epst} that the resulting sequence is still compact in $L^2$.  Thus we can take a limit and obtain a nontrivial solution $\eps$ to the \emph{linearised} equation
$$ \eps_t + \eps_{xxx} + (p R^{p-1} \eps)_x = \alpha(t) R_x $$
for some scalar quantity $\alpha(t)$, which stays compact in $L^2$ and bounded in $H^1$ (and obeys the orthogonality conditions \eqref{rtx}).  The task is thus to show that there is no such solution other than the trivial solution $\eps = 0$, which will establish the Liouville theorem and thus Theorem \ref{asymt}.

At this point, one can now use the global virial estimate \eqref{epsor}, which is valid here due to the exponential decay of $\eps$.  If $\eps$ is non-trivial, it has an $H^1$ norm bounded away from zero, which in conjunction with the $L^2$ compactness shows that the right-hand side of \eqref{epsor} is negative and bounded away from zero.  But this forces $\int_\R (x-x(t)) \eps^2$ to go to $-\infty$, which contradicts the exponential decay of $\eps$.  This finally finishes the argument.  (An alternate proof of this ``linear Liouville theorem'' was also recently established in \cite{martel-linear}.)

We conclude with a number of miscellaneous remarks.

\begin{remark}
This proof was significantly more complicated than the direct proof, but the underlying strategy is much more powerful: it uses compactness methods to strip away all the inessential portions of the dynamics, leaving a very smooth and localised solution to which global estimates can be applied.  This dispenses with the need for any cutoffs in space or frequency, which can significantly complicate the analysis.
\end{remark}

\begin{remark}
These arguments have recently been simplified and generalised further (to handle arbitrary power-type nonlinearities) in \cite{mm-general}. Interestingly, this result does not require any spectral hypotheses on the linearised problem.  Also, some sharper asymptotics for special nonlinearities (in particular the $p=4$ case) have been obtained in \cite{tao:gkdv4}, \cite{martel2008}.
\end{remark}

\begin{remark} The above arguments crucially rely on the tendency of gKdV evolutions to propagate non-soliton mass in one direction, while the soliton itself moves in another direction, leading to some crucial monotonicity formulae.  These methods can be extended to some other models \cite{dika}, \cite{dika2}, \cite{mizu}.  For NLS models however, there are only weak analogues of these formulae (see \cite{mmt2}), although for some special nonlinearities (in which spectral hypotheses on the linearised operator are assumed) asymptotic stability for such equations can be recovered \cite{bus}, \cite{perelman}, \cite{rss}.  The situation is particularly well understood in the radial case, in which the soliton remains tied to the spatial origin and so the degeneracies in stability associated with translation invariance are eliminated.  For a survey of these issues see \cite{cuccagna}.
\end{remark}

\begin{remark}
In a sort of converse to asymptotic stability (analogous to the existence of wave operators in scattering theory), it is also possible to construct solutions with prescribed asymptotic behaviour of the type described above; see \cite{cote}, \cite{cote2}, \cite{cote3}.
\end{remark}

\section{The critical case}\label{crit}

We now discuss the more difficult mass-critical case $p=5$, in which the scale invariance now plays a much more delicate role. The details are rather technical and we shall only paraphrase or sketch the key arguments here, referring the reader to the original papers for full details.
For further surveys of the results in this section, see \cite{mm-survey}, \cite{mm-survey2}, \cite{tzvetkov}.

The new difficulty in the critical case can be seen by considering the one-parameter family of soliton states $Q_c(x) := c^{1/(p-1)} Q(c^{1/2} x)$.  A short computation shows that $M(Q_c) = c^{\frac{2}{p-1}-\frac{1}{2}} M(Q)$ and $E(Q_c) = c^{\frac{2}{p-1}+\frac{1}{2}} E(Q)$.  If $p \neq 5$, the exponents here are non-zero (and $M(Q)$ and $E(Q)$ are also non-zero), and so the conservation of mass or energy prohibits a solution which starts near $Q$ from ending up near $Q_c$ for any $c \neq 1$.  But in the critical case $p=5$, we have $M(Q_c) = M(Q)$, while one can compute that $E(Q_c)=E(Q)=0$.  Thus, the mass and energy conservation laws do not prohibit the possibility of drift from $Q$ to $Q_c$, or more generally
along the entire ground state surface \eqref{sigmap}, on which the mass is always $M(Q)$ and the energy is always $0$.  Indeed, numerics \cite{bona-numerics} suggest that solutions $u$ starting near a ground state will increase their scale parameter $c$ to infinity in finite time, thus leading to finite time blowup (thus, for instance $\lim_{t \to T^-} \|u(t)\|_{H^1_x} = +\infty$ for some finite time $T$).  An inspection of the linearised operator also supports the possibility of drift to increasingly finer scales.

Henceforth we fix $p=5$.  As discussed earlier, the sharp Gagliardo-Nirenberg inequality of Weinstein\cite{Weinstein} allows one to show that no blowup occurs as long as $M(u) < M(Q)$, in particular there are initial data arbitrarily close to the ground state for which one has global existence.  A refinement of this analysis also allows one to consider the situation\footnote{The case $M(u)=M(Q)$ was subsequently treated in \cite{mm-duke}, in which blowup was shown to be impossible.} in which $M(Q) < M(u) < M(Q)+\alpha$ for some small $\alpha > 0$, as long as the energy $E(u)$ is \emph{negative} (this situation can again occur arbitrarily close to the ground state $Q$, as long as the mass is strictly greater than $M(Q)$).  In this case, we can again (up to a harmless change of sign, $u \mapsto -u$) obtain a decomposition of the form \eqref{utx-asym} over the lifespan of the solution.  However, a key difference in the critical case is that the scale parameter $c(t)$ can go to infinity in finite time.

Henceforth we assume $M(Q) < M(u) < M(Q)+\alpha$ for sufficiently small $\alpha$, and also $E(u) < 0$. It will be convenient to make the error $\eps$ \eqref{utx-asym} dimensionless (i.e. invariant under the scaling symmetry) by replacing \eqref{utx-asym} with the equivalent decomposition
$$
u(t,x) = c(t)^{1/4} [ Q( c(t)^{1/2}( x - x(t) ) )  + \eps( t, c(t)^{1/2}( x - x(t) ) ) ].$$
As before, one can select two orthogonality conditions on $\eps$; it turns out to be convenient to require $\langle \eps(t), Q^3 \rangle = \langle \eps(t), Q' \rangle = 0$.  In that case one can use Weinstein's analysis to show (after replacing $u$ with $-u$ if necessary) that $\eps(t)$ is small in $H^1$ (see e.g. \cite{merle} for details\footnote{In the paper \cite{merle} and the other papers cited here, the wavelength parameter $\lambda(t) := c(t)^{-1/2}$ is used rather than the scale parameter $c(t)$, but this of course makes only a minor notational difference to the argument.}).  In particular we see that $\|u(t)\|_{H^1_x} \sim c(t)^{1/2}$, and so blowup in the $H^1$ norm is equivalent to $c(t)$ going to infinity.

Suppose for the moment that the solution $u$ existed globally in time, with $c(t)$ bounded both above and below, and furthermore that $\eps(t)$ ranged in a compact subset of $L^2$ (so in particular, one has bounds such as \eqref{xl}, \eqref{xr}).  Then the same Liouville theorem analysis used in the subcritical case can be used to show that $\eps = 0$, which would of course contradict the assumption that $M(u)$ is strictly greater than $M(Q)$; see \cite{liouville} for details.  By repeating the rest of the subcritical analysis, this is already enough to deduce asymptotic stability in the case when $c(t)$ is bounded both above and below.  However, this statement is vacuous due to the results in \cite{merle}, which in fact show that $c(t) \to \infty$ as $t \to T$, where $0 < T \leq \infty$ is the maximal time of existence.

We sketch the proof of this result as follows.  Once again we use a compactness and contradiction argument.  If the above claim failed, then one could find a sequence $u = u_n$ of solutions, each with some finite or infinite lifespan $0 < T_n \leq +\infty$, with mass $M(Q) < M(u_n) < M(Q)+o(1)$ and $E(u_n)<0$, such that the velocity function $c_n(t)$ stayed bounded in $t$ for each $n$ (though with a bound depending on $n$).  For each $n$, we consider the quantity $c_{n,*} := \liminf_{t \to T_n} c_n(t)$.  Then $c_{n,*}$ is bounded; from the negative conserved energy one can also show $c_{n,*}$ to be non-zero.  Thus we can find a time $t_n$ for each $n$ such that $c_n(t_n)$ is very close to $c_{n,*}$, and that $c_n(t)$ is either close to or larger than $c_{n,*}$ for all $t > t_n$.  By rescaling and translating in time if necessary we may take $c_n(t_n)=1$ and $t_n=0$.

For each $n$, we see from construction that $c_n(t)$ is always greater than or close to $1$, and returns for some infinite sequence of times $t_{n,m} \to T_n$ to be close to $1$.  This situation is similar to the previous situation ``$c_n(t)$ bounded above and below'', except that we now allow $c_n$ to oscillate between being close to $1$ and being extremely large.  It turns out that with enough care, one can use local mass propagation estimates much as before to obtain exponential decay similar to \eqref{epst} at and near the times $t_{n,m}$, or more precisely for a limiting error profile $\eps_{n,\infty}$ obtained as a weak limit of $\eps_n(t+t_{n,m},x)$.  In particular, this places the limiting error profile $\eps_{n,\infty}$ (and thus the limiting solution profile $u_{n,\infty}$) in $L^1$.  This allows one to deploy a final conservation law, namely the mean $\int_\R u_{n,\infty}(t)$.  The conservation of this quantity can be used to show that the limiting velocity parameter $c_\infty(t)$ is bounded both above and below, at which point the Liouville theorem ensures tht $\eps_{n,\infty}$ must vanish.  But this turns out to be incompatible with the negative energy hypothesis.  This contradiction establishes the desired claim that one has blowup either at finite or at infinite times.

It turns out that one can analyse these solutions further.  For this it turns out to be convenient to change the orthogonality conditions slightly, so that one now requires
$$ \langle \eps(t), \frac{Q}{2} + yQ_y \rangle = \langle \eps(t), y(\frac{Q}{2} + yQ_y) \rangle = 0$$
where $y$ denotes the spatial variable.  (The expression $\frac{Q}{2} + yQ_y$ can be viewed as the infinitesimal scaling vector field at $Q$.)  These conditions are convenient for applying virial identities.  It is also convenient to introduce a dimensionless time variable $s$ defined by the ODE \begin{equation}\label{dsdt}
\frac{ds}{dt} = c^{3/2}(t),
\end{equation}
One can show that regardless of whether $T$ is finite or infinite, that $s \to +\infty$ as $t \to T$.  The reason for these choices of coordinates is that the error $\eps = \eps(s,y)$ now obeys an analytically tractable PDE, which takes the form
$$ \eps_s = (L \eps)_y + N(\eps,Q)$$
where $L$ is the linear operator $L := -\partial_{xx} + 1 - 5Q^4$ and $N$ is a nonlinear (and slightly non-local) expression which consists of terms which are quadratic and higher in $\eps$, and which are largely localised in space (due to the fact that most terms involve at least one factor of the localised function $Q$).  The finer structure of $N$ is important for certain computations, but for simplicity we will not delve into the explicit form of $N$ here.

If $c$ is now viewed as a function of $s$ rather than $t$, then from the preceding analysis we have $c(s) \to \infty$ as $s \to \infty$.  However, this divergence to infinity could potentially be quite oscillatory; there is no \emph{a priori} reason why $c(s)$ should be monotone.  Also we do not have \emph{a priori} knowledge as to the rate of growth of $c$ in $s$.

Suppose however that we had a time interval $[s_1,s_2]$ (which could potentially be very long) in which $c(s)$ varies between $c(s_1)$ and $c(s_2) = 1.1 c(s_1)$ (say), thus there is a slight focusing effect along this interval.  Suppose also that $c(s) > c(s_1)$ for all times $s \geq s_2$, and also that we have an exponential localisation such as \eqref{epst} on the time interval $[s_1,s_2]$ (actually, for this argument, we only need exponential localisation to the left of the soliton).  It turns out that there is no solution of this form (assuming our hypotheses of negative energy, and mass close to that of the ground state, of course).  This key fact is established in \cite{mm-annals}, and we sketch the proof as follows.  We may rescale so that $c(s_1)=1$, thus $c(s) \sim 1$ for $s \in [s_1,s_2]$.  The localisation places $\eps$ in $L^1$ (to the left, at least).  We introduce a quantity $J$ that measures the amount of $L^1$ mass of $\eps$ that is to the left of the soliton; the precise definition of $J$ is
\begin{equation}\label{jdef}
 J(s) := \int \eps(s,y) (\int_y^{+\infty} (\frac{Q}{2} + z Q_z)\ dz)\ dy - \frac{1}{4} (\int Q)^2.
\end{equation}
It turns out that $J$ obeys a differential equation of the form
$$ \frac{d}{ds} (c^{-1/4} J) = - 2 c^{-1/4} \int \eps Q + O( c^{-1/4} \int \eps^2 e^{-|y|^2/2}\ dy ),$$
which reflects the fundamental fact that the non-soliton portion of mass tends to propagate to the left. 
On the other hand, the conservation of mass and energy eventually gives us a lower bound
$$ \int \eps Q \geq \frac{1}{c} |E(u)| + \frac{1}{4} \int \eps_y^2\ dy + O( \int \eps^2 e^{-|y|^2/2}\ dy )$$
and so on the interval $[s_1,s_2]$ (where $c$ lies between $1$ and $1.1$) we have
\begin{equation}\label{diffeq}
\frac{d}{ds} (c^{-1/4} J) + 0.1 \int \eps_y^2\ dy + 0.1 |E(u)| \leq O( \int \eps^2 e^{-|y|^2/2}\ dy )
\end{equation}
(say).  On the other hand, on $[s_1,s_2]$ the quantity $c$ increases from $1$ to $1.1$, which together with \eqref{jdef} tells us that $c^{-1/4} J$ increases over this interval by at least some absolute constant $\delta > 0$ (recall from \eqref{epst} that $\eps$ is small in $L^1$).  Thus if we integrate \eqref{diffeq} over $[s_1,s_2]$ we arrive at a bound of the form
$$ 1 + \int_{s_1}^{s_2} \int \eps_y^2\ dy + |E(u)| |s_2-s_1| \leq O( \int_{s_1}^{s_2} \int \eps^2 e^{-|y|^2/2}\ dy ).$$
On the other hand, virial identity arguments (similar to those used to derive \eqref{epsx}) combined with the hypothesis that $c(s)$ is comparable to $1$ on $[s_1,s_2]$, allow us to establish a bound of the form
$$ \int_{s_1}^{s_2} \int \eps^2 e^{-|y|^2/2}\ dy \leq C \alpha^{1/2} ( 1 + \int_{s_1}^{s_2} \int \eps_y^2\ dy + |E(u)| |s_2-s_1| )$$
where $\alpha := M(u) -M(Q)$.  For $\alpha$ sufficiently small, this leads to the desired contradiction.

By combining the above result with the Liouville theorem discussed earlier, and another compactness-and-contradiction argument, it was shown in \cite{mm-annals} that there existed a sequence of times $s_m \to \infty$ such that $\eps(s_m)$ converged weakly (in $H^1$) to zero.  Indeed, since we know $c(s) \to \infty$ as $s \to \infty$ and is continuous in $s$, we can define $s_m$ to be the last time for which $c(s_m) = (1.1)^m$.  If $\eps(s_m)$ did not converge weakly to zero, then (after passing to a subsequence if necessary) it would converge to some other limiting error profile $\eps_\infty(0)$, which can be viewed as the error at time zero for some limit solution profile $u_\infty(s)$, with its attendant velocity parameter $c_\infty(s)$, which is bounded from below by $c_\infty(0)$ for all $s \geq 0$ by construction.  A careful application of local mass propagation inequalities (see \cite{mm-annals} for details) allows one to establish exponential decay on $\eps_\infty(s)$ on the left.  If $c_\infty(s)$ increases to $1.1 c_\infty(0)$ in finite time, one can then use the previous argument to obtain a contradiction.  If instead $c_\infty(s)$ ranged between $c_\infty(0)$ and $1.1 c_\infty(0)$ for all $s$, then one can use further local mass propagation to get some decay on the right as well.  The Liouville theorem then would force $\eps_\infty$ to vanish, again leading to a contradiction.  Thus we obtain the weak convergence of $\eps(s_m)$ to $0$ as claimed.

A further compactness-and-contradiction argument in \cite{mm-annals} allows one to strengthen this claim, to assert that $\eps(s)$ converges weakly to $0$ for \emph{all} times $s \to \infty$, not just a subsequence $s_m \to \infty$ of times.  For if this were not the case, one could find another sequence $s'_m \to \infty$ of times for which $\eps(s'_m)$ was converging to a non-zero limit $\eps'_\infty(0)$, associated to a limiting solution profile $u'_\infty$.  But because any sufficiently large $s'_m$ can be placed between a pair $s_{m'}, s_{m'+1}$ of times where $\eps$ is close to zero (in the weak topology), it is possible to use local mass conservation to show that $M(u'_\infty) \leq M(Q)$; meanwhile, the energy $E(u'_\infty)$ can be shown by limiting arguments to be non-positive.  But this, together with the orthogonality conditions on $\eps'_\infty$, forces $\eps'_\infty$ to vanish, a contradiction.

Finally, in \cite{mm-jams}, it was shown that solutions of the above form in fact blow up in finite time (thus giving theoretical confirmation of the numerical blowups observed for instance in \cite{bona-numerics}), provided that one make an additional assumption of decay on the right (e.g. it would suffice to have $\int_0^\infty u(0,x)^2 x^{6+\delta}\ dx < \infty$ for some $\delta > 0$).  

We very briefly sketch the main ideas as follows.  In view of the relation \eqref{dsdt}, some numerology shows that finite time blowup will eventuate if we can show that $\frac{d}{ds} \log c(s)$ grows (on average at least) like $c(s)^{-3/2+\delta}$ for some $\delta > 0$.  In fact, in \cite{mm-jams} it is shown that
\begin{equation}\label{cs}
\frac{d}{ds} \log c(s) \geq \delta |E(u)| c(s)^{-1}
\end{equation}
``on the average'', for some absolute constant $\delta > 0$, which leads to finite time blowup for negative energy data.

A calculation of the dynamics of $c(s)$ reveals an equation of the form
$$ \frac{d}{ds} \log c(s) = 2 \int \eps L((\frac{Q}{2} + y Q_y)_y)\ dy + O( \int \eps^2 e^{-|y|/100}\ dy ).$$
The error term here is controllable by virial identities.  But the quantity $\int \eps L((\frac{Q}{2} + y Q_y)_y)\ dy$ is somewhat oscillatory and is not easy to control directly.  To avoid this problem, the arguments in \cite{mm-jams} introduce a new decomposition
$$
u(t,x) = c'(t)^{1/4} [ Q( c'(t)^{1/2}( x - x'(t) ) )  + \eps'( t, c'(t)^{1/2}( x - x'(t) ) ) ]$$
where $\eps'$ obeys some rather different orthogonality conditions, namely
$$ \int (\int_{-\infty}^y \frac{Q}{2} + Q_z\ dz) \eps'(y)\ dy = \int y (\frac{Q}{2} + Q_y) \eps'\ dy = 0.$$
Because $\int_{-\infty}^\infty \frac{Q}{2} + Q_z\ dz$ is non-zero, such orthogonality conditions are only reasonable when $u$ has sufficient decay on the right.  Fortunately, the hypothesis of polynomial decay on the right at time zero, together with some local mass propagation estimates, turn out to give enough decay to make this decomposition viable (see \cite{mm-jams} for details).  We also have an associated rescaled time variable $s'$ defined by $ds'/dt = c'(t)^{3/2}$.  The point of performing this decomposition is that the equation for $c'(s')$ is simpler than that for $c(s)$:
$$ \frac{d}{ds'} \log c'(s') = 2 \int \eps' Q + O( \int (\eps')^2 e^{-|y|/100}\ dy ).$$
Furthermore, one can show that $c$ and $c'$ are comparable, as are $\eps$ and $\eps'$, in various technical senses.  In particular, the quantities $\int \eps' Q$ and $\int \eps Q$ can be related to each other modulo acceptable errors.  On the other hand, mass and energy conservation considerations eventually let one establish a lower bound on $\int \eps Q$ of the form $\delta |E(u)|/c(s)$, modulo acceptable error terms.  The claim follows.

\begin{remark}
The above analysis in fact gives some more explicit upper and lower bounds on the blowup rate (i.e. the behaviour of $c(t)$ as $t$ approaches $T$); see \cite{mm-jams} for details.
\end{remark}

\begin{remark}
For the analogous problem for mass-critical NLS, blowup from negative energy data can be obtained easily from a virial identity argument, at least when the data is localised in space; see \cite{glassey-blow}, while explicit blowup solutions can also be constructed via the pseudo-conformal symmetry of this equation.  However, the corresponding virial identities for gKdV equation, while leading to important estimates such as \eqref{epsx}, seem to be unable to coerce blowup just by themselves.
\end{remark}

\begin{remark}
The situation for the supercritical equation $p>5$ is still poorly understood; due to some crucial changes of sign, many of the methods used above fail completely.  Numerics such as those in \cite{bona-numerics} continue to suggest finite time blowup starting near soliton initial data in this case, but the dynamics are known to be unstable \cite{bss} and so it is unlikely that a controlled blowup of the type seen above will hold, at least for generic data.
\end{remark}

\section{Further developments}

We now briefly survey more recent developments regarding stability of multisolitons, and on collisions between solitons.  Due to the vast explosion in activity in these areas (for a wide variety of dispersive models), we will not be able to give a comprehensive bibliography here, instead focusing on a representative subset of results.

There has been some progress in understanding the stability aspects for \emph{multisolitons} - superpositions of two or more solitons.  As the underlying equations are non-linear, constructing these solutions is non-trivial, and indeed existence of such solutions is usually established simultaneously with stability results.  When the solitons are far apart and receding from each other, one is in a perturbative regime and can in some cases glue together the stability theory for single solitons to create multisolitons.  See \cite{mmt} for an instance of this for subcritical gKdV multisolitons, and \cite{rss}, \cite{mmt2}, \cite{perelman} for some results for subcritical NLS multisolitons.

More difficult is the question of what happens when two solitons collide.  Recently there has been some work on the collision between fast thin solitons and slow broad solitons \cite{martel-collision}, \cite{martel-collision2}; the situation is still perturbative, but there are noticeable nonlinear effects, such as a shift in position in the fast soliton caused by nonlinear interaction with the slow soliton.  Somewhat similarly in spirit, there has been some recent work in \cite{hz}, \cite{holmer} investigating the collision between a fast soliton and a localised potential (such as a delta function potential).  The fully non-perturbative situation of two large slow solitons colliding is however still beyond current technology.

In supercritical cases, solitons are generally unstable.  However in some cases the number of unstable directions is finite, and so finite dimensional stable manifolds can (in principle) be constructed.  This turns out to be a rather delicate issue (relying, among other things, on strong control on the spectrum of the linearised operator, which can contain some non-trivial resonances); see \cite{krieger}, \cite{krieger-stable} for recent developments.

There has also been many recent papers on the blowup in the neighbourhood of a soliton to a scale-invariant evolution equation; these results do not quite follow exactly the same pattern as the results for critical gKdV mentioned in the previous section, but certainly share many of the same ingredients.  For instance, see \cite{mr0}, \cite{mr}, \cite{mr2} for some results relating to the mass-critical NLS.  For the energy-critical non-linear wave equation or wave maps equation, there are some slightly different approaches to create blowup \cite{krieger2}, \cite{kst}, \cite{rs}; see \cite{schlag} for a recent survey.

The compactness-and-contradiction approach has also been applied to derive the \emph{absence} of blowup (i.e. global existence) for critical equations in the \emph{defocusing} case, or in focusing cases in which the solution is ``smaller'' in mass or energy than that of the ground state.  See \cite{merlekenig}, 
\cite{tvz-higher}, \cite{KTV}

\bibliographystyle{amsplain}

\end{document}